\documentclass{article}
\usepackage{fullpage,amssymb}
\def\proof{{\em Proof:\/} \hspace{0.3cm}}
\newtheorem{theorem}{Theorem}

\begin{document}

\title{On Indecomposable Normal Matrices \\
in Spaces with Indefinite Scalar Product}
\author{O.V. Holtz \\
\small Department of Applied Mathematics \\
\small Chelyabinsk State Technical University \\
\small 454080 Chelyabinsk, Russia}
\date{}
\maketitle

\begin{abstract}
 Finite dimensional linear spaces (both complex and real) with
indefinite scalar product $[\cdot, \cdot]$ are considered.  Upper and
lower bounds are given for the size of an indecomposable matrix that
is normal with respect to this scalar product in terms of specific
functions of $v=\min\{v_{-},v_{+}\}$, where $v_{-}$ ($v_{+}$) is the
number of negative (positive) squares of the form $[x,x]$. All the
bounds except for one are proved to be strict.
\end{abstract}

\section{Definitions and notation}

Consider a complex (real) linear space $C^{n}$ ($R^{n}$) with an
indefinite scalar product $[ \cdot , \cdot ]$. By definition, the
latter is a nondegenerate sesquilinear (bilinear) Hermitian form.  If
the usual scalar product $(\cdot , \cdot)$ is fixed, then there exists
a nonsingular Hermitian operator $H$ such that $[x,y]=(Hx,y)$ $\forall
x,y\in C^{n}$ ($R^{n}$).  If $A$ is a linear operator, then the
{\em{$H$-adjoint\/}} of $A$ (denoted by $A^{[*]}$) is defined by the
identity $[A^{[*]}x,y] \equiv [x,Ay]$. An operator $N$ is called
{\em{$H$-normal\/}} if $NN^{[*]}=N^{[*]}N$. An operator $U$ is called
{\em{$H$-unitary\/}} if $UU^{[*]}=I$, where $I$ is the identity
transformation.

Let $V$ be a nontrivial subspace of $C^{n}$ ($R^{n}$). The subspace
$V$ is called {\em{neutral}} if $[x,y]=0$ $\forall x,y \in V$.  If the
conditions $x \in V$ and $[x,y]=0$ $\forall y \in V$ imply $x=0$, then
$V$ is called {\em{nondegenerate}}. The subspace $V^{[\perp]}$ is
defined as the set of all vectors $x$ from $C^{n}$ ($R^{n}$) such that
$[x,y]=0$ $\forall y \in V$. If $V$ is nondegenerate, then
$V^{[\perp]}$ is also nondegenerate and $V \dot{+} V^{[\perp]}=C^{n}$
($R^{n}$), where $\dot{+}$ stands for the direct sum.

A linear operator $A$ is called {\em{decomposable}} if there exists
a nondegenerate proper subspace $V$ of $C^{n}$ ($R^{n}$) such that
both $V$ and $V^{[\perp]}$ are invariant under $A$ or (it is the
same) if $V$ is invariant both under $A$ and $A^{[*]}$. Then $A$
is the $H$-{\em{orthogonal sum of $A_{1}=A|_{V}$ and
$A_{2}=A|_{V^{[\perp]}}$}}. If an operator $A$ is not decomposable,
it is called {\em{indecomposable}}.

By the {\em{rank}} of a space we mean $v=\min\{v_{-},v_{+}\}$,
where $v_-$ ($v_+$) is the number of negative (positive)
squares of the form $[x,x]$, i.e., the number
of negative (positive) eigenvalues of the operator $H$.

The problem is to find  functions $f_1(\cdot)$, $f_2(\cdot)$
such that $f_1(v) \leq n \leq f_2(v)$ for any indecomposable
$H$-normal operator acting in a space of dimension $n$ and of rank
$v$ and to find out whether these bounds are strict.

This problem arises in the classification of indecomposable
$H$-normal matrices \cite{2,3}. The bounds for the size of an
indecomposable $H$-normal matrix in a complex space are known
\cite{2}. In Section~2, we check their strictness. The bounds for
matrices in real spaces are considered in Section~3.

As in \cite{2} and \cite{3}, we denote by $I_{r}$ the $r \times r$
identity matrix, by $D_{r}$ the $r \times r$ matrix with 1's on the
trailing diagonal and zeros elsewhere, and by
$A \oplus B \oplus \ldots \oplus Z$ the block diagonal matrix with
blocks $A$, $B$, $\ldots$, $Z$. By $A^T$ we mean $A$ transposed.

\section{Indecomposable normal matrices in complex spaces}
The objective of this section is to prove the following theorem.
\begin{theorem}
Let an indecomposable $H$-normal operator $N$ act in a space $C^{n}$
of rank $k>0$. Then either (A) or (B) holds:

\begin{description}
\item{(A)} $N$ has only one eigenvalue and $2k \leq n \leq 4k$;
\item{(B)} $N$ has only two eigenvalues and $n=2k$;
\end{description}
these bounds being strict.
\end{theorem}

\proof Theorem~1 of \cite{2} states that for an indecomposable
$H$-normal operator $N$ there exist two alternatives: (A) and (B) so
that it suffices to prove that these estimates are unimprovable.

{\em{Step 1.\/}} Show the strictness of the low bound in (A), i.e.,
for any $k>0$ point out a pair of $2k \times 2k$ matrices
$\{N,H\}$, where $H$ has $k$ negative and $k$ positive eigenvalues,
$N$ is $H$-normal and indecomposable and has only one eigenvalue
$\lambda$. Let
\begin{equation}
N=\left(\begin{array}{cc} \lambda I_{k} & N_{1} \\ 0 & \lambda I_{k}
      \end{array} \right), \;\;\;\;\;\; H=\left(\begin{array}{cc} 0 &
      I_{k} \\ I_{k} & 0 \end{array} \right).  \label{NH-2k-1}
\end{equation}
It can easily be checked that $N$ is $H$-normal.
In addition, let the submatrix $N_{1}$ be nonsingular.

Proposition~1 from \cite{2} and its Corollary may be restated as
follows:

{\em{Let an $H$-normal operator $N$ acting in $C^{n}$ have $\lambda$
as its only eigenvalue. Let the subspace \begin{equation} S_{0}=\{ x
\in C^{n}: \; (N - \lambda I)x=(N^{[*]}- \overline{\lambda} I)x=0\}
\label{S0} \end{equation} be neutral. Then there exists a
decomposition of $C^{n}$ into a direct sum of subspaces $S_{0}$, $S$,
$S_{1}$ such that
\begin{equation}
N=\left( \begin{array}{ccc} N'=\lambda I & * & * \\ 0 & N_{1} & * \\ 0
        & 0 & N''=\lambda I \end{array} \right), \;\;\;\;\;\; H=\left
        ( \begin{array}{ccc} 0 & 0 & I \\ 0 & H_{1} & 0 \\ I & 0 & 0
        \end{array} \right), \label{pred1} \end{equation} where $N':$
        $S_{0} \rightarrow S_{0}$, $N_{1}:$ $S \rightarrow S$, $N'':$
        $S_{1} \rightarrow S_{1}$, the internal operator $N_{1}$ is
        $H_{1}$-normal, and the pair $\{N_{1},H_{1}\}$ is determined
        up to the unitary similarity. To go over from one
        decomposition $C^{n}=S_{0}\dot{+}S\dot{+}S_{1}$ to another by
        means of a transformation $T$ it is necessary that $T$ be
        block triangular with respect to both decompositions.\/}}

\noindent (In Proposition~1 from~\cite{2} there are two conditions:

(a) $N$ is indecomposable

(b) $n>1$

\noindent instead of the condition

(c) $S_{0}$ is neutral,

\noindent but the results are only derived from (c), which follows
from (a) and (b)).

We see that (\ref{NH-2k-1}) is a specific case of (\ref{pred1})
corresponding to the decomposition $C^{n}=S_{0}\dot{+}S\dot{+}S_{1}$
with $S=0$. If there exists a nondegenerate subspace $V$ such that
both $V_{1}=V$ and $V_{2}=V^{[\perp]}$ are invariant under $N$,
then, according to our restatement of Proposition~1 from \cite{2},
for $i=1,2$ we have
$V^{(i)}=S_{0}^{(i)}\dot{+}S^{(i)}\dot{+}S_{1}^{(i)}$,
where $S_{0}^{(i)}=V^{(i)} \cap S_{0}$ and the pairs
$\{N^{(i)}, H^{(i)}\}$ have the form (\ref{pred1}). But
(\ref{NH-2k-1}) implies $S_{0}=S_{0}^{[\perp]}$ so that for any
$i=1,2$ the subspace $S^{(i)}$ is trivial. Thus, $N$ from
(\ref{NH-2k-1}) is decomposable if and only if there exists a
transformation $T$ preserving $H$ and reducing $N$ to the form
$$  \widetilde{N}=\left( \begin{array}{cc}
        \lambda I_{k} & \widetilde{N_{1}} \\
        0 & \lambda I_{k}
      \end{array}
   \right), \]
where $\widetilde{N_{1}}$ is block diagonal (that is, $T$ is an
$H$-unitary transformation of $N$ to the form $\widetilde{N}$).
The matrix $T$ is necessarily block triangular with respect to the
decomposition $C^{n}=S_{0}\dot{+}S_{1}$, i.e.,
$$  T=\left( \begin{array}{cc}
        T_{1} & T_{2} \\
        0 & T_{3}
      \end{array}
   \right). \]
 For $T$ to be $H$-unitary it is necessary to have
$T_{3}=T_{1}^{*-1}$. Then from the condition $NT=T\widetilde{N}$
it follows that $N_{1}=T_{1}\widetilde{N_{1}}T_{1}^{*}$.
Therefore, $N$ will be indecomposable if $N_{1}$ is not
congruent to any block diagonal matrix $\widetilde{N_{1}}$.

If $N_{1}$ and a block diagonal matrix $\widetilde{N_{1}}$ are
congruent, then $N_{1}N_{1}^{*-1}$ is similar to
$\widetilde{N_{1}}\widetilde{N_{1}}^{*-1}$.
Since the latter is also block diagonal, for $N_{1}$ to be not
congruent to $\widetilde{N_{1}}$ it is sufficient that
$N_{1}N_{1}^{*-1}$ cannot be reduced to block diagonal form.

Let us prove that for any $n=1,2, \ldots$ there exists a nonsingular
real $(2n-1)\times(2n-1)$ matrix $N_{1}$ such that
\begin{equation}
 N_{1}N_{1}^{*-1}=\left( \begin{array}{cccccc}
        1 & 1 & 0 & \ldots & 0 & 0 \\
        0 & 1 & 1 & \ldots & 0 & 0 \\
        \vdots & \vdots & \vdots & \ddots & \vdots & \vdots \\
        0 & 0 & 0 & \ldots & 1 & 1 \\
        0 & 0 & 0 & \ldots & 0 & 1 \end{array}
   \right)  \label{N1-odd}
\end{equation}
and a nonsingular real $2n \times 2n$ matrix $N_{1}$ such that
\begin{equation}
  N_{1}N_{1}^{*-1}=\left( \begin{array}{cccccc}
        -1 & 1 & 0 & \ldots & 0 & 0 \\
        0 & -1 & 1 & \ldots & 0 & 0 \\
        \vdots & \vdots & \vdots & \ddots & \vdots & \vdots \\
        0 & 0 & 0 & \ldots & -1 & 1 \\
        0 & 0 & 0 & \ldots & 0 & -1 \end{array}
   \right)  \label{N1-even}
\end{equation}

(now it is not necessary that $N_{1}$ be real, but this will be
used in the next section). The matrices (\ref{N1-odd}) and
(\ref{N1-even}) are obviously not similar to any block diagonal
ones because for each of them the subspace generated by their
eigenvectors is one-dimensional.

Prove the statement by induction for odd numbers.  If $n=1$, let
$N_{1}^{(1)}=(1)$. Suppose we have found a nonsingular real
$(2n-1)\times(2n-1)$ matrix $N_{1}^{(n)}$ with the property
required. Let
$$  N_{1}^{(n+1)}=\left( \begin{array}{ccc}
        0 & A_{n+1} & B_{n+1} \\ C_{n+1} & N_{1}^{(n)} & 0 \\ D_{n+1}
        & 0 & 0 \end{array} \right), \] where the submatrices
        $A_{n+1}$, $B_{n+1}$, $C_{n+1}$, $D_{n+1}$ will be shortly
        specified.  If we denote by $\Lambda_{n}$ the $(2n-1) \times
        (2n-1)$ matrix (\ref{N1-odd}), by $\Lambda_{n}'$ the
        $1\times(2n-1)$ matrix
$$  \Lambda_{n}'=\left( \begin{array}{ccccc}
        1 & 0 & 0 & \cdots & 0 \end{array}
   \right), \]
and by $\Lambda_{n}''$ the $(2n-1) \times 1$ matrix
$$  \Lambda_{n}''=\left( \begin{array}{ccccc}
        0 & \cdots & 0  & 0 & 1 \end{array}  \right)^T, \]
then the condition $N_{1}^{(n+1)}=\Lambda_{n+1}N_{1}^{(n+1)*}$
may be rewritten as follows:
\begin{eqnarray}
0 & = & \Lambda_{n}'A_{n+1}^{*}, \label{cong1} \\
A_{n+1} & = & C_{n+1}^{*}+\Lambda_{n}'N_{1}^{(n)*}, \label{cong2} \\
B_{n+1} & = & D_{n+1}^{*}, \label{cong3} \\
C_{n+1} & = & \Lambda_{n} A_{n+1}^{*}+ \Lambda_{n}''B_{n+1}^{*}
\label{cong4}
\end{eqnarray}
(by the inductive hypothesis,
$N_{1}^{(n)}=\Lambda_{n}N_{1}^{(n)*}$). Taking
\begin{eqnarray*}
& A_{n+1}=\left( \begin{array}{ccccc}
        0 & a_{2} & a_{3} & \ldots & a_{2n-1} \end{array}
   \right), \;\;\;\;\; B_{n+1}=(b), & \\
& a_2, \; a_3, \; \ldots, \; a_{2n-1}, \; b \in \Re, &
\end{eqnarray*}
$C_{n+1}=\Lambda_{n}A_{n+1}^{*}+\Lambda_{n}''B_{n+1}^{*}$,
$D_{n+1}=B_{n+1}^{*}$, one can satisfy the conditions (\ref{cong1}),
(\ref{cong3}), (\ref{cong4}). Substituting the expression for
$C_{n+1}$ in (\ref{cong2}), we get the only condition to be
satisfied:
\begin{equation}
-N_{1}^{(n)}\Lambda_{n}'^{*}=(\Lambda_{n}-I)A_{n+1}^{*}+
\Lambda_{n}''B_{n+1}^{*}.  \nonumber
\end{equation}
Since its right hand side is equal to
$$  \left( \begin{array}{ccccc}
a_{2} &   a_{3} & \cdots & a_{2n-1} &  b \end{array}  \right)^T, \]
it always can be satisfied by choosing the appropriate values of
$a_{2}$, $a_{3}$, $\ldots$, $a_{2n-1}$, $b$. By construction, the
elements of each matrix $N_{1}^{(n)}$ below the trailing diagonal
are zeros and those on the trailing diagonal are equal to $\pm 1$
so that $N_{1}^{(n)}$ is nonsingular. Thus, for odd numbers our
statement is proved.

To prove it for even numbers one can take
$$ N_{1}^{(1)}= \left( \begin{array}{cc}
       \frac{1}{2} & 1 \\
       -1 & 0 \end{array}  \right) \]
and construct $N_{1}^{(n+1)}$ from $N_{1}^{(n)}$ in the same way
as before (the details are left to the reader). Step~1 is completed.

{\em{Step 2.}} Show the strictness of the upper bound in (A).
The example of the pair $\{N,H\}$ is
\begin{equation}
N=\left( \begin{array}{cccc}
        \lambda I_{k} & I_{k} & 0 & 0 \\
        0 & \lambda I_{k} & 0 & N_{1} \\
        0 & 0 & \lambda I_{k} & N_{2} \\
        0 & 0 & 0 & \lambda I_{k} \end{array}
   \right), \;\;\;\;
 H=\left( \begin{array}{cccc}
        0 & 0 & 0 & I_{k} \\
        0 & I_{k} & 0 & 0 \\
        0 & 0 & I_{k} & 0 \\
        I_{k} & 0 & 0 & 0 \end{array}
   \right),  \label{NH-4k}
\end{equation}
where
$$ N_{1}=\left( \begin{array}{ccccc}
        0 & r_{1} & 0 & \ldots & 0 \\
        0 & 0 & r_{2} & \ldots & 0 \\
        \vdots & \vdots & \vdots & \ddots & \vdots \\
        0 & 0 & 0 & \ldots & r_{k-1} \\
        r_{k} & 0 & 0 & \ldots & 0 \end{array}
   \right), \]
$$ N_{2}=\left( \begin{array}{cccc}
        \sqrt{1-r_{1}^{2}} & 0 & \ldots & 0 \\ 0 & \sqrt{1-r_{2}^{2}}
        & \ldots & 0 \\ \vdots & \vdots & \ddots & \vdots \\ 0 & 0 &
        \ldots & \sqrt{1-r_{k}^{2}} \end{array} \right), \] $r_{i} \in
        (0,1)$ $\forall i=1, \ldots, k$, and $r_{i} \neq r_{j}$ if
        $i\neq j$. The matrix $H$ has $k$ negative and $3k$ positive
        eigenvalues. The matrix $N$ is $H$-normal, since the condition
        $N_{1}^{*}N_{1}+N_{2}^{*}N_{2}=I$ is satisfied. As before, we
        see that (\ref{NH-4k}) is a specific case of
        (\ref{pred1}). Suppose a nondegenerate subspace $V$ is
        invariant both under $N$ and under $N^{[*]}$. Denote the basis
        vectors of $S_{0}$ by $\{v_{i}\}_{i=1}^{k}$, the basis vectors
        of $S_{1}$ by $\{w_{i}\}_{i=1}^{k}$ (here the basis
        corresponds to (\ref{NH-4k})). Let $\widetilde{V}$ be the
        range of the projection of $V$ onto $S_1$ along
        $S_0^{[\perp]}$. It is a subspace of dimension $m>0$, since
        $V$ necessarily contains at least one nontrivial vector from
        $S_{0}$ and, therefore, at least one vector with nontrivial
        projection onto $S_{1}$ (otherwise $V$ would be
        degenerate). Let
        $\{\sum_{j=1}^{k}\alpha_{ij}w_{j}\}_{i=1}^{m}$ be a basis of
        $\widetilde{V}$. From the condition $(N-\lambda
        I)(N^{[*]}-\overline{\lambda}I)V \subseteq V$ it follows that
        $\{\sum_{j=1}^{k}\alpha_{ij}v_{j}\}_{i=1}^{m} \subset V$.  If
        $V$ is nondegenerate, $\widetilde{V}$ and
        $\widetilde{S_{0}}=S_{0} \cap V$ are necessarily of the same
        dimension. Therefore,
        $\{\sum_{j=1}^{k}\alpha_{ij}v_{j}\}_{i=1}^{m}$ is a basis of
        $\widetilde{S_{0}}$. As $(N -\lambda I)^{2}V \subseteq V$,
        $(N^{[*]} -\overline{\lambda} I)^{2}V \subseteq V$, we obtain
        $\{\sum_{j=1}^{k}\alpha_{ij}N_{1}v_{j}\}_{i=1}^{m} \subset V$,
        $\{\sum_{j=1}^{k}\alpha_{ij}N_{1}^{*}v_{j}\}_{i=1}^{m} \subset
        V$.  As $V^{[\perp]} \cap S_{0} \neq \{0\}$, we have $m(=dim
        \widetilde{S_{0}})<k$. Thus, for $N$ to be decomposable it is
        necessary that the subspace $\widetilde{S_{0}}$, which is of
        dimension more than zero and less than $k$, be invariant under
        $N_{1}$ and under $N_{1}^{*}$. This means the existence of an
        orthogonal projection $P (\neq 0, \; I)$ commuting with
        $N_{1}$.  But it can easily be checked by direct calculation
        that from the conditions $N_{1}P=PN_{1}$ and $P=P^{*}$ it
        follows that $P=\mu I$.  Since $P^{2}=P$, we have $\mu =0$ or
        $\mu=1$ so that $P=0$ or $P=I$.  The contradiction obtained
        shows that $N$ is indecomposable. Step~2 is completed.

{\em{Step 3.}} Now for any $k>0$ let us point out a pair of
$2k \times 2k$ matrices $\{N,H\}$, where $H$ has $k$ negative
and $k$ positive eigenvalues, $N$ is $H$-normal and indecomposable
and has the two eigenvalues $\lambda_{1}$, $\lambda_{2}$
($\lambda_{1} \neq \lambda_{2}$). Let
\begin{equation}
N=\left( \begin{array}{cc}
        N_{1} & 0 \\
        0 & N_{2} \end{array}
   \right), \;\;\;\;\;\;\;
 H=\left( \begin{array}{cc}
        0 & I_{k} \\
        I_{k} & 0 \end{array}
   \right),  \label{NH-2k-2}
\end{equation}
where the $k \times k$ matrices $N_{i}$ ($i=1,2$)
are as follows:
$$ N_{1}=\left( \begin{array}{ccccc}
        \lambda_{1} & 1 & 0 & \ldots & 0 \\
        0 & \lambda_{1} & 1 & \ldots & 0 \\
        \vdots & \vdots & \vdots & \ddots & 1 \\
        0 & 0 & 0 & \ldots & \lambda_{1} \end{array}
   \right), \;\;\;\; N_{2}=\lambda_{2} I_{k}. \]
The matrix $N$ is $H$-normal, for it satisfies the condition
$N_{1}N_{2}^{*}=N_{2}^{*}N_{1}$.
Suppose that $N$ is similar to
$$ \widetilde{N}=\left( \begin{array}{cc}
        N|_{V}=\widetilde{N_{1}} & 0  \\
        0 & N|_{V^{[\perp]}}=\widetilde{N_{2}} \end{array}
   \right), \]
where $V$ is a nondegenerate subspace. Since the subspace generated
by the eigenvectors of $N$ corresponding to the eigenvalue
$\lambda_{1}$ is one-dimensional, one of the submatrices
$\widetilde{N_{1}}$, $\widetilde{N_{2}}$ (for example,
$\widetilde{N_{2}}$) has $\lambda_{2}$ as its only eigenvalue,
hence $\widetilde{N_{2}}=\lambda_{2} I$. But any subspace
generated by eigenvectors corresponding to the eigenvalue
$\lambda_{2}$ is neutral so that $V$ cannot be nondegenerate. This
contradiction shows the indecomposability of $N$. Step~3 is
completed. The theorem is proved. \hfill $\Box$

\section{Indecomposable normal matrices in real spaces}
The objective of this section is to prove the following theorem.
\begin{theorem}
Let an indecomposable $H$-normal  operator $N$ act in a space $R^{n}$
of rank $k>0$. Then one of the conditions (A) - (E) holds:
\begin{description}
\item{(A)} $N$ has only one real eigenvalue and $2k \leq n \leq 4k$;
\item{(B)} $N$ has only two real eigenvalues and $n=2k$;
\item{(C)} $N$ has only two complex conjugate eigenvalues and
       $n=2$ if $k=1$ and   $2k \leq n \leq 10[k/2]-2$ if $k>1$;
\item{(D)} $N$ has only one real and one pair of complex conjugate
eigenvalues and $n=2k$;
\item{(E)} $N$ has only two pairs of complex conjugate eigenvalues
and $n=2k$.
\end{description}
The alternatives (D) and (E) are possible only if $k$ is even. The
estimates (A), (B), (D), (E), and the low bound in (C) are strict.
\end{theorem}

\proof That an indecomposable $H$-normal matrix has one of
the five sets of eigenvalues is proved in~\cite[Lemma~1]{3}.
Bounds (A) and (B) are proved in~\cite[Theorem~1]{2},
their strictness in Theorem~1 from the previous section (since the
matrices constructed in Theorem~1 are real and any matrix that is
indecomposable in a complex space is also indecomposable in a real
one). The condition $n \geq 2k$ is obvious. Indeed, since
$k=\min\{v_{-},v_{+}\}$ and $n=v_{-}+v_{+}$, we have $n \geq 2k$.
Thus, we must consider the cases (C) - (E) only, keeping in mind
that $n \geq 2k$.

{\em{Step 1.}} Consider the case (C). Let $N$ have the two distinct
eigenvalues $\lambda=\alpha+i\beta$ and
$\overline{\lambda}=\alpha-i\beta$. The equality $n=2$ for $k=1$ is
proved in~\cite[Theorem~1]{3}. In case when $k=2$ Theorem~2
of~\cite{3} states that $n \leq 8$. So, it remains to prove the
inequality $n \leq 10[k/2]-2$ for $k \geq 3$. To this end recall
Proposition~2 from~\cite{3}:

{\em{ Let an indecomposable $H$-normal operator $N$ acting in
$R^{n}$ ($n>2$) have the two distinct eigenvalues
$\lambda=\alpha+i\beta$, $\overline{\lambda}=\alpha-i\beta$. Let
$$S_{0}'=\{z=x+iy \; (x,y \in R^{n}): \;  Nz=\lambda z,
\; N^{[*]}z=\overline{\lambda}z \}, \]
$$S_{0}''=\{z=x+iy \; (x,y \in R^{n}): \;  Nz=\lambda z,
\; N^{[*]}z=\lambda z \}, \]
$\{z_{j}\}_{1}^{p}$ ($\{z_{j}\}_{p+1}^{p+q}$) be a basis of
$S_{0}'$ ($S_{0}''$), and
\begin{equation}
S_{0}=\sum_{j=1}^{p+q}span\{x_{j},y_{j}\}. \label{S0-2ev}
\end{equation}
Then there exists a decomposition of $R^{n}$ into a direct sum of
subspaces $S_{0}$, $S$, $S_{1}$ such that
\begin{equation}
 N=\left( \begin{array}{ccc}
 N' & * & *    \\
 0 & N_{1} & * \\
 0 & 0 & N''
 \end{array} \right), \;\;\;\;\;\;
   H=\left( \begin{array}{ccc}
 0 & 0 & I    \\
 0 & H_{1} &  0    \\
 I & 0 & 0 \end{array} \right), \label{aug_dec1}
\end{equation}
where
$$ N': \; S_{0} \rightarrow S_{0}, \;\;\;\;\;\;
N'=N_{1}' \oplus \ldots \oplus N_{p+q}', \]
\begin{equation}
 N_{j}'=\left( \begin{array}{cc}
 \alpha & \beta    \\
 -\beta & \alpha \end{array} \right), \;\;\;\;\;\;  j=1, \ldots p+q,
 \label{aug_dec2}
\end{equation}
$$ N'': \; S_{1} \rightarrow S_{1}, \;\;\;\;\;\;
N''=N_{1}'' \oplus \ldots \oplus N_{p+q}'', \]
\begin{equation}
 N_{j}''=N_{j}' \;\; if \;1 \leq j \leq p, \;\;
 N_{j}''=N_{j}'^{*} \;\; if \; p < j \leq p+q, \label{aug_dec3}
\end{equation}
the internal operator $N_{1}$ is $H_{1}$-normal and the pair
$\{N_{1},H_{1}\}$ is determined up to unitary similarity.
To go over from one decomposition $R^{n}=S_{0}\dot{+}S\dot{+}S_{1}$
to another by means of a transformation $T$ it is necessary that
the matrix $T$ be block triangular with respect to both
decompositions.}}

According to this proposition, for an indecomposable
operator $N$ the subspace $S_{0}$ defined in (\ref{S0-2ev})
is neutral so that its dimension does not exceed $k$.
Therefore, if we prove that for $n>10[k/2]-2$ the condition
$dim\:S_{0} \leq k$ fails, this will mean the decomposability
of $N$.

According to~\cite[the proof of Lemma~1]{1}, if an $H$-normal
operator $N$ acting in $C^{n}$ has the two distinct eigenvalues
$\lambda$, $\overline{\lambda}$, then there exists a decomposition
of $C^{n}$ into a direct sum of subspaces $V_{1}$, $V_{2}$,
$V_{3}$, $V_{4}$ such that
$$ N=\left( \begin{array}{cccc}
          N|_{V_{1}}=N_{1} & 0 & 0 & 0  \\
          0 & N|_{V_{2}}=N_{2} & 0 & 0  \\
          0 & 0 & N|_{V_{3}}=N_{3} & 0  \\
          0 & 0 & 0 & N|_{V_{4}}=N_{4}  \end{array} \right), \]
$$ H=\left( \begin{array}{cccc}
          0 & I & 0 & 0 \\
          I & 0 & 0 & 0 \\
          0 & 0 & H_{3} & 0 \\
          0 & 0 & 0 & H_{4} \end{array} \right). \]
Here $N_{1}$ and $N_{3}$ have only one eigenvalue $\lambda$, $N_{2}$
and $N_{4}$ only one eigenvalue $\overline{\lambda}$, and
$dim\:V_{1}=dim\:V_{2}$. In our case $C^{n}$ is $R^{n}$ complexified,
therefore, $dim\:V_{3}=dim \: V_{4}$ too. Either $V_{1}$ or $V_{3}$
may be equal to zero.

Let $n > 10[k/2]-2$, i.e, $n \geq 10[k/2]$. Consider the
following three cases: (a) $V_{1}=V_{2}=0$, (b) $V_{3}=V_{4}=0$,
(c) $dim \:V_{1}>0$ and $dim V_{3}>0$.

(a) If $V_{1}=0$, then $dim \:V_{3}(=dim\:V_{4}) \geq 5[k/2]$.
Let $H_3$ ($H_4$) have $v_{-(3)}$ ($v_{-(4)}$) negative eigenvalues.
Without loss of generality it can be assumed that
$k=v_-=v_{-(3)}+v_{-(4)}$ so that
$\min \{v_{-(3)},v_{-(4)}\} \leq [k/2]$. Let $v_{-(3)} \leq [k/2]$.
Decompose $N_{3}$ into an $H$-orthogonal sum of indecomposable
operators $N_{3}^{(1)}$, $N_{3}^{(2)}$, $\ldots$, $N_{3}^{(m)}$:
$N_{3}=N_{3}^{(1)}\oplus N_{3}^{(2)}\oplus \ldots \oplus
N_{3}^{(m)}$, $H_{3}=H_{3}^{(1)}\oplus H_{3}^{(2)}\oplus
\ldots \oplus H_{3}^{(m)}$, $V_{3}=V_{3}^{(1)}\dot{+}
V_{3}^{(2)}\dot{+} \ldots \dot{+} V_{3}^{(m)}$.
Denote by $v_{-(3)}^{(j)}$ the number of negative eigenvalues of
 $H_{3}^{(j)}$ ($j=1, \ldots, m$). Let
$$ V_{3}'=\sum_{v_{-(3)}^{(j)}>0} V_{3}^{(j)}, \;\;\;\;\;\;
 V_{3}''=\sum_{v_{-(3)}^{(j)}=0} V_{3}^{(j)}, \]
$H_{3}'$, $H_{3}''$ and $V_{3}'$, $V_{3}''$ be the corresponding
sums of $H_{3}^{(j)}$ and $V_{3}^{(j)}$. Since for
$v_{-(3)}^{(j)}>0$ the condition $dim V_{3}^{(j)} \leq 4
v_{-(3)}^{(j)}$ holds \cite[Theorem~1]{2},
we have $dim \:V_{3}' \leq 4v_{-(3)} \leq 4[k/2]$.
If $v_{-(3)}^{(j)}=0$, then $dim\:V_3^{(j)}=1$,
$H_3^{(j)}=(1)$, and $N_{3}^{(j)}=(\lambda)$.
Thus, $dim V_{3}'' \geq [k/2]$ and $Nz=\lambda z$,
$N^{[*]}z=\overline{\lambda}z$ for all $z \in V_{3}''$.
The operators $N_{3}'$ and $N_{3}'^{[*]}$ commute so that
if $dim\:V_{3}'>0$, there exists at least one vector
$z_{0} \in V_{3}'$ such that $Nz_{0}=\lambda z_{0}$ and
$N^{[*]}z_{0}=\overline{\lambda}z_{0}$. If $dim \:V_{3}'=0$, all
nontrivial vectors from $V_{3}$ are eigenvectors of $N$
corresponding to the eigenvalue $\lambda$ and those of $N^{[*]}$
corresponding to the eigenvalue $\overline{\lambda}$. Therefore,
in either case there exist at least $p=[k/2]+1$ linearly
independent vectors $\{z_{l}\}_{l=1}^{p}$ such that
$Nz_{l}=\lambda z_{l}$, $N^{[*]}z_{l}=\overline{\lambda}z_{l}$.
Therefore, $dim S_{0}\geq 2([k/2]+1)>k$.

(b) If $dim V_{3}=0$, then $n=2\:dim\:V_{1}$. Since no neutral
subspace of a space of rank $k$ can be of dimension more than
$k$, $dim\:V_{1} \leq k$ so that $n \leq 2k$. But it was proved
before that $n \geq 2k$. Thus, in this case $n=2k < 10[k/2]-2$.

(c) If $dim V_{1}, \;dim V_{3}>0$, we can assume, as in the case
(a) above, that $v_{-(3)} \leq [k/2]-dim V_1$ (the notation here
is also as in (a)). Since $dim \:V_{3} \geq 5[k/2]-dim\:V_{1}$,
there are $p$ linearly independent vectors $\{z_{l}\}_{l=1}^{p}$
such that $Nz_{l}=\lambda z_{l}$,
$N^{[*]}z_{l}=\overline{\lambda}z_{l}$
($p$ is equal to $[k/2]+3 dim V_1 +1$ if $dim V_{3}'>0$ or to
$5[k/2]-dim V_{1}$ if $dim V_{3}'=0$).
Since $dim V_{1} \leq k$ and $k \geq 3$, we have
$5[k/2]-dim V_{1} \geq [k/2] +1$.
Thus, again $dim S_{0} >k$ so that $N$ is decomposable.
The upper bound in (C) is proved. Step 1 is completed.

{\em{Step 2.}} Let us show the strictness of the low bound in
(C) for even numbers $k$. Consider the pair of $2k \times 2k$
matrices
$$ N=\left( \begin{array}{cccccc}
          A & I_{2} & 0 & \ldots & 0 & 0  \\
          0 & A & I_{2} & \ldots & 0 & 0  \\
          \vdots & \vdots & \vdots & \ddots & \vdots & \vdots  \\
          0 & 0 & 0 & \ldots & A & I_{2}  \\
          0 & 0 & 0 & \ldots & 0 & A  \end{array} \right), \;\;\;\;
 H=\left( \begin{array}{ccccc}
          0 & 0 & \ldots & 0 & I_{2}  \\
          0 & 0 & \ldots & I_{2} & 0  \\
          \vdots & \vdots & \ddots & \vdots & \vdots  \\
          0 & I_{2} & \ldots & 0 & 0  \\
          I_{2} & 0 & \ldots & 0 & 0  \end{array} \right), \]
where
\begin{equation}
 A=\left( \begin{array}{cc}
          \alpha & \beta  \\
          -\beta & \alpha \end{array} \right) \;\;\;\;\;
(\alpha, \beta \in \Re, \; \beta > 0) \label{A}
\end{equation}
(throughout what follows, by $A$ we will denote the matrix
(\ref{A})). It is seen that $H$ has $k$ negative and $k$ positive
eigenvalues. It can easily be checked by direct calculation that
$N$ is $H$-normal. The number of linearly independent vectors
$z_{l}$ satisfying the condition $Nz_{l}=\lambda z_{l}$
($\lambda = \alpha+i\beta$) is equal to $1$, hence $dim S_{0}=2$.
By~\cite[Proposition~3]{3}, if the subspace $S_{0}$ is
two-dimensional, the operator $N$ is indecomposable. So, the
statement is proved and Step~2 is completed.

{\em{Step 3.}} For the case when $k$ is odd consider the following
pair of $2k \times 2k$  matrices
$$  N=\left( \begin{array}{ccccccccc}
          A & X & \cdots & 0 & 0 & 0 & \cdots & 0 & 0 \\
          0 & A & \cdots & 0 & 0 & 0 & \cdots & 0 & 0 \\
  \vdots & \vdots & \cdot & \vdots & \vdots & \vdots &
  \cdot & \vdots & \vdots \\
          0 & 0 & \cdots & A & X & 0 & \cdots & 0 & 0 \\
          0 & 0 & \cdots & 0 & A & X & \cdots & 0 & 0 \\
          0 & 0 & \cdots & 0 & 0 & A^{*} & \cdots & 0 & 0 \\
  \vdots & \vdots & \cdot & \vdots & \vdots & \vdots &
  \cdot & \vdots & \vdots \\
          0 & 0 & \cdots & 0 & 0 & 0 & \cdots & A^{*} & X \\
          0 & 0 & \cdots & 0 & 0 & 0 & \cdots & 0 & A^{*}
\end{array} \right), \]
$$  H=\left( \begin{array}{ccccccccc}
          0 & 0 & \cdots & 0 & 0 & 0 & \cdots & 0 & I_2 \\
          0 & 0 & \cdots & 0 & 0 & 0 & \cdots & I_2 & 0 \\
  \vdots & \vdots & \cdot & \vdots & \vdots & \vdots &
  \cdot & \vdots & \vdots \\
          0 & 0 & \cdots & 0 & 0 & I_2 & \cdots & 0 & 0 \\
          0 & 0 & \cdots & 0 & D_2 & 0 & \cdots & 0 & 0 \\
          0 & 0 & \cdots & I_2 & 0 & 0 & \cdots & 0 & 0 \\
  \vdots & \vdots & \cdot & \vdots & \vdots & \vdots &
  \cdot & \vdots & \vdots \\
          0 & I_2 & \cdots & 0 & 0 & 0 & \cdots & 0 & 0 \\
          I_2 & 0 & \cdots & 0 & 0 & 0 & \cdots & 0 & 0
\end{array} \right), \]
where
$$ X=\left( \begin{array}{cc}
          1 & 1 \\
          1 & 1 \end{array} \right). \]
As $N^{[*]}=N$, the matrix $N$ is $H$-normal. Since the condition
$Nz=\lambda z$ ($\lambda = \alpha + i \beta$) implies
$$ z=\left( \begin{array}{ccccc}
 z_1 & iz_1 & 0 & \cdots & 0
\end{array} \right)^{T},$$ 
the subspace $S_0$ is two-dimensional and, according
to~\cite[Proposition~3]{3}, the matrix $N$ is indecomposable.
Step~3 is completed.

{\em{Step 4.}} Consider the case (D). Let $N$ have one real
eigenvalue $\lambda$ and two complex conjugate eigenvalues
$\alpha \pm i \beta$. According to~\cite[Proposition~1]{3},
$R^{n}$ is a direct sum of neutral subspaces ${\cal Q}_{1}$,
${\cal Q}_{2}$ such that $dim {\cal Q}_{1}=dim {\cal Q}_{2}$,
$N{\cal Q}_{1} \subseteq {\cal Q}_{1}$,
$N{\cal Q}_{2} \subseteq {\cal Q}_{2}$,
$N|_{{\cal Q}_{1}}$ has $\lambda$ as its only eigenvalue,
$N|_{{\cal Q}_{2}}$ has the two eigenvalues $\alpha \pm i \beta$.
 From the last condition it follows that $dim {\cal Q}_{2}$
is even. As in Step 1, case (b), we have
$n=2\:dim{\cal Q}_{1} \leq 2k$, hence $n=2k$ and
$dim {\cal Q}_{2}=k$ is necessarily an even number.

Now suppose that k is even and consider the following
pair $\{N,H\}$:
$$ N=\left( \begin{array}{cc}
        N_{1} & 0 \\
        0 & N_{2} \end{array}
   \right),  \;\;\;\;\;\;
H=\left( \begin{array}{cc}
        0 & I_{k} \\
        I_{k} & 0 \end{array}
   \right),  \]
where  the $k\times k$ submatrices $N_{1}$ and $N_{2}$ are as
follows:
$$ N_{1}=\left( \begin{array}{cccccc}
        A & I_{2} & 0 & \ldots & 0 & 0 \\
        0 & A & I_{2} & \ldots & 0 & 0 \\
        \vdots & \vdots & \vdots & \ddots & \vdots & \vdots \\
        0 & 0 & 0 & \ldots & A & I_{2} \\
        0 & 0 & 0 & \ldots & 0 & A \end{array}   \right),
\;\;\;\;\;\; N_{2}=\lambda I_{k}.  \]
It is clear that the condition
\begin{equation}
N_{1}N_{2}^{*}=N_{2}^{*}N_{1} \label{norm-last}
\end{equation}
is satisfied so that $N$ is $H$-normal. Suppose that there exists a
nondegenerate subspace $V$ such that $N$ is similar to the matrix
\begin{equation}
 \widetilde{N}=\left( \begin{array}{cc}
        N|_{V}=\widetilde{N_{1}} & 0  \\
        0 & N|_{V^{[\perp]}}=\widetilde{N_{2}} \end{array}
   \right).  \label{n-new}
\end{equation}
Since the subspace generated by the eigenvectors corresponding to
the eigenvalue $\alpha + i \beta$ is one-dimensional (in the
complexified space), one of the submatrices $\widetilde{N_{1}}$
and $\widetilde{N_{2}}$ has $\lambda$ as its only eigenvalue,
therefore, either $\widetilde{N_{1}}$ or $\widetilde{N_{2}}$ is
equal to $\lambda I$. As in Theorem~1, Step~3, we conclude that
under this condition $V$ cannot be nondegenerate so that $N$ is
indecomposable. Step 4 is completed.

{\em{Step 5.}} Consider the case (E). That $k$ is necessarily even
and $n=2k$ can be proved just as in Step 4 before.
So, it remains to construct a pair $\{N, H\}$ satisfying the
$H$-normality condition, where $H$ has $k$ negative and $k$ positive
eigenvalues, $N$ is indecomposable, and $N$ has only two pairs of
complex conjugate eigenvalues $\alpha_{1} \pm i\beta_{1}$,
$\alpha_{2} \pm i\beta_{2}$ ($\beta_{1}, \; \beta_{2}>0$,
$(\alpha_{1}, \beta_{1}) \neq (\alpha_{2}, \beta_{2})$).
Let
$$ N=\left( \begin{array}{cc}
        N_{1} & 0 \\
        0 & N_{2} \end{array}
   \right),  \;\;\;\;\;\;
H=\left( \begin{array}{cc}
        0 & I_{k} \\
        I_{k} & 0 \end{array}
   \right),  \]
where  the $k\times k$ submatrices $N_{1}$ and $N_{2}$ are as
follows:
$$ N_1=\left( \begin{array}{cccccc}
        A_1 & I_2 & 0 & \ldots & 0 & 0 \\
        0 & A_1 & I_2 & \ldots & 0 & 0 \\
        \vdots & \vdots & \vdots & \ddots & \vdots & \vdots \\
        0 & 0 & 0 & \ldots & A_1 & I_2 \\
        0 & 0 & 0 & \ldots & 0 & A_1 \end{array}   \right), \;\;
 N_2=\left( \begin{array}{cccc}
        A_2 & 0 & \ldots & 0 \\
        0 & A_2 & \ldots & 0 \\
        \vdots & \vdots & \ddots & \vdots \\
        0 & 0 & \ldots & A_2 \end{array} \right).  \]
Here
$$ A_1=\left( \begin{array}{cc}
        \alpha_1 & \beta_1 \\
        -\beta_1 & \alpha_1 \end{array}
   \right),  \;\;\;\;\;\;
A_2=\left( \begin{array}{cc}
        \alpha_2 & \beta_2 \\
        -\beta_2 & \alpha_2 \end{array}
   \right).  $$
It can easily be checked that the $H$-normality condition
(\ref{norm-last}) is satisfied. As in Step~4 before, the
assumption that $N$ is similar to (\ref{n-new}) implies that either
$\widetilde{N_{1}}$ or $\widetilde{N_{2}}$ (suppose
$\widetilde{N_{2}}$) has two eigenvalues $\alpha_{2} \pm i\beta_{2}$
only. Therefore, there are $m=\frac{dim \widetilde{N_{2}}}{2}$
complex linearly independent eigenvectors
$\{z_{j}=x_{j}+iy_{j}\}_{j=1}^{m}$
of $\widetilde{N_{2}}$ corresponding
to the eigenvalue $\alpha_{2} + i\beta_{2}$. Consequently, the
set $\{x_{j}\}_{j=1}^{m} \cup \{y_{j}\}_{j=1}^{m}$ is a basis of
$\widetilde{N_{2}}$. But
$[x_{j}, x_{l}]=[y_{j},y_{l}]=[x_{j},y_{l}]=0$
for all $j,l=1,\ldots m$. Therefore, the subspace $V$ cannot be
nondegenerate and hence $N$ is indecomposable. Step~5 is completed.
The theorem is proved. \hfill $\Box$

{\em{I would like to thank Prof. Vladimir Strauss for constant
attention to my work and the referee for valuable suggestions.}}

\end{document}